\documentclass[12pt]{article}
\usepackage{amsmath,amsfonts,amsthm,amscd,upref,amstext}
\usepackage[dvips]{graphicx}
\usepackage{subfigure}

\newtheorem{teore}{Theorem}[section]
 \newtheorem{lema}[teore]{Lemma}
 \newtheorem{propo}[teore]{Proposition}
 \newtheorem{coro}[teore]{Corollary}
\newtheorem{defi}[teore]{Definition}
\newtheorem{rem}[teore]{Remark}

\newtheorem*{agra}{Acknowlegment}

\newcommand{\R}{\mathbb R}
\newcommand{\C}{\mathbb C}
\newcommand{\cala}{\mathcal A}
 \newcommand{\cali}{\mathcal I}
 \newcommand{\calo}{\mathcal O}
\newcommand{\calr}{\mathcal R}

\begin{document}

\title{Topological Triviality of Families of Singular Surfaces}

\author{R. Callejas-Bedregal\thanks{Work partially supported by
 PADCT/CT-INFRA/CNPq/MCT - Grant 620120/04-5}\,,\,\,K. Houston\thanks{The second
 named author thanks the ICMC at S\~ao Carlos for their
hospitality. His research was partially supported by EPSRC, Grant
reference GR/S48639/01}\,\,\,and \,M. A. S. Ruas\thanks{The third
named author thanks the ICTP for its hospitality and financial
support during her visit to Trieste, and CNPq, Brazil, grant
\#300066/88-0 and FAPESP, grant \#03/03107-9, for partial support.
2000 Mathematics Subject Classification: 58K40(primary),
58K65(secondary).}}

\date{}
\maketitle

\section{Introduction}
\label{sec:intro}  We study the topological triviality of families
of singular surfaces in $\C^3$ parame\-trized by $\mathcal
A$-finitely determined map germs.

Finitely determined map germs $f:(\mathbb \C^{2},0)\rightarrow
(\mathbb \C^{3},0)$ can be approximated by stable maps $f_{s},$
and information on the topology of such approximations can be
obtained in terms of data calculable from the original map germ
$f$.   In \cite{mond}, D.~Mond defines the $0$-stable invariants,
$T(f)$, the number of triple points of $f_{s}$, and $C(f)$, the
number of Whitney umbrellas of $f_{s}$,  and shows how to compute
them in terms of $f$. These two invariants together with
$\mu(D^{2}(f))$, the Milnor number of the double point locus, form
a complete set of invariants for  $\cala$-simple  germs
(\cite{mond}).

The following natural extension of the above result was formulated
by Mond in \cite{mond1}. Let $f_{t} : (\C^{2},0)\rightarrow
(\C^{3},0)$ be a one-parameter family of finitely determined map
germs. Does the constancy of the invariants $T(f_{t})$, $C(f_{t})$
and $\mu({D}^{2}(f_{t}))$ imply the topological triviality of the
family? In this work we give a positive answer to this question as
a consequence of the following L\^e-Ramanujam type theorem,
ensuring that the constancy of the Milnor number of the double
point locus characterizes the topological triviality of the
family:

\begin{teore}\label{topological}
Let $F:(\C^2\times\C,0) \to (\C^3\times\C,0)$ be a one-parameter
unfolding of an  $\cala$-finitely determined map germ $f: (\C^2,0)
\to (\C^3,0).$ The following statements are equivalent:

\begin{enumerate}
\item  $\mu(D^{2}(f_{t}))$ is constant for $t\in T,$ where $T$ is
a small neighborhood of $0$ in $\C;$ \item  $F$ is topologically
trivial.
\end{enumerate}
\end{teore}

 We  prove that if the unfolding  $F$ is $\mu$-constant, in the
sense that $\mu(D^{2}(f_{t}))$ is constant for $t\in T,$ where $T$
is a small neighborhood of $0$ in $\C,$ then $F$ is excellent, as
defined by T. Gaffney  in \cite{Gaffney}. An excellent unfolding
is a Thom stratified mapping and we obtain the trivialization by
integrating controlled vector fields tangent to the strata of the
stratification of $F$ given by the stable types in source and
target.

A natural question is whether $\mu$-constant in a one-parameter
family $F$ also implies the Whitney equisingularity of the family.
This is indeed the case, and it follows as a consequence of the
following theorem which completely describes the equisingularity
of $F$ in terms of the equisingularity of $D^2(F),$ the double
point locus in source, that happens to be a family of reduced
plane curves. There is a well-understood theory of equisingularity
for plane curves, starting with the results of O. Zariski
(\cite{Za}), further developed among others, by Zariski himself,
H.~Hironaka, M.~Lejeune, L\^e D\~ung Tr\`ang and B.~Teissier (see
\cite{GaffMass,Tei} for surveys on the subject and \cite{FernBob}
for some interesting new developments). For families of reduced
plane curves, it is well known that topological triviality,
Whitney equisingularity and bilipschitz equisingularity are
equivalent notions. It is very surprising that this is  also true
for families of singular surfaces in $\C^3$ parametrized by
$\cala$-finitely determined map germs, as shown by the following
theorem:

\begin{teore}\label{equisingularities}
Let $F:(\C^2\times\C,0) \to (\C^3\times\C,0)$ be a one-parameter
unfolding of an  $\cala$-finitely determined map germ $f: (\C^2,0)
\to (\C^3,0).$ Then, the following are equivalent:

\begin{enumerate}
\item  $F$ is topologically trivial. \item $F$ is Whitney
equisingular. \item $F$ is bilipschitz trivial.

\end{enumerate}

\end{teore}

Theorems \ref{topological} and \ref{equisingularities} will be
proved in Theorems \ref{topological1} and \ref{bilipschitz},
respectively.

For other results related to the subject discussed in this paper,
see for instance \cite{jonu,Ruas}.

\begin{agra}{\rm
The authors thank A.C.G.\ Fernandes, V.H.\ Jorge-P\'erez and J.J.\
Nu\~no-Ballesteros for helpful conversations. They especially
thank L\^e Dung Tr\'ang for suggesting the usage of vector fields
in the proof of the main result.}
\end{agra}

\section{Previous results}

 We first review some  results on the geometry and
classification of singularities of surfaces in 3-space.

\begin{defi}
  Two map germs $f, g:(\C^n,0)\to (\C^p,0)$ are $\cala-$equivalent,
  denoted by $g\sim_\cala f,$ if
  there exist map germs of diffeomorphisms $h:(\C^n,0)\to (\C^n,0)$
   and $k:(\C^p,0)\to (\C^p,0),$
   such that $g=k\circ f \circ h^{-1}.$
  \end{defi}

  \begin{defi}
$f:(\C^n,0)\to (\C^p,0)$ is finitely determined ($\cala$-finitely
determined) if there exists a positive integer $k$ such that for
any $g$ with $j^kg(0)=j^kf(0)$ we have $g\sim_\cala f.$
\end{defi}


\begin{teore}[(Geometric criterion, Mather-Gaffney, 1975)] A map germ
$f:(\C^n,0)\to (\C^p,0)$ is finitely determined if, and only if,
for every representative $f$ (of $f$) there exist $U$, a
neighborhood of $0$ in $\C^n$, and  $V$, a neighborhood of $0$ in
$\C^p,$ with $f(U)\subset V$, such that for all $y\in
V\setminus{0},$ the set $S=f^{-1}(y)\cap \Sigma(f)$ is finite and
$f:(\C^n,S) \to (\C^p,y)$ is stable.
\end{teore}

Let $f:(\C^{2},0)\rightarrow (\C^{3},0)$ be a finitely determined
map germ. From the classical result of  Whitney \cite{W}, we know
that the stable singularities in these dimensions are transverse
double points, triple points and cross-caps. In this case, the
above theorem implies that $f:(\C^2,0)\to (\C^3,0)$ is
$\cala$-finitely determined if and only if for every
representative $f$ there exists a neighborhood $U$ of $0,$ such
that the only singularities of $f(U)\setminus \{0\}$ are
transverse double points.

A finitely determined map germ  $f:(\C^2,0)\to (\C^3,0)$ has a
versal unfolding

$$\begin{array}{ccl}
  F:(\C^2 \times\C^r,0) & \longrightarrow  & (\C^3\times\C^r,0) \\
 (x,u)                  & \mapsto          & (f_u(x),u).
\end{array}$$

 Given a representative of the map germ $F$, which we also denote by
 $ F,$ defined in a neighborhood $U \times W$ of $(0,0)$ in
 $\C^2\times \C^r,$ we can define the bifurcation set
$B\subset \C^r$  of $F,$ by $B=\{u\in W: f_u\ \text{is not
stable}\}.$

The set $B$ is a proper algebraic subset of $\C^r,$  hence its
complement is a connected set in $\C^r.$ Hence, for any $u,u' \in
W \setminus B,$ $f_u$ and $f_{u'}$ are stable and $f_u\sim_\cala
f_{u'}.$

 Notice that $f_u$ has a finite number of cross-caps and
triple points and these are analytic invariants of the original
map germ.

Let $f_u$ be a local stable perturbation of $f$, (known as the
{\em{disentanglement}} of $f$).
 In \cite{mond} D.~Mond defines the following $0$-stable invariants
of $f$:
\begin{eqnarray*}
C(f)&=& {\mbox{\# of cross-caps of }} f_u , \\
T(f)&=& {\mbox{\# of triple-points of }} f_u .
\end{eqnarray*}


Formulas to compute $C(f)$ and $T(f)$ as the codimension of
certain algebras associated to $f$ are given in \cite{mond}.





\section{Double point locus}

To study the topology of $f(\C^2),$ in a small neighborhood of
$0,$ one needs a third invariant associated to  the double point
locus of $f$ which we now describe. Let $f:U\to \C^p$ be a
holomorphic map, where $U\subset \C^n$ is an open subset and
$n\leq p.$ We define the {\it double point set of} $f$, denoted by
$\tilde{D}^{2}(f),$ as the closure in $U\times U$ of the set
$$\{(x, y)\in U\times U: f(x)=f(y), x\neq y\} .$$

To choose a convenient analytic structure for the double point set
$\tilde{D}^{2}(f),$ we follow the construction of
\cite{Mond-Marar} which is also valid for holomorphic maps from
$\C^n$ to $\C^p$ with $n\leq p.$ Let us denote the diagonals in
$\C^n\times \C^n$ and $\C^p\times \C^p$ by $\Delta_n$ and
$\Delta_p$ respectively and denote the sheaves of ideals defining
them by ${\cali}_n$ and ${\cali}_p$ respectively. We write the
points of $\C^n\times \C^n$ as $(x, x').$ Then, for each $i=1,
\ldots, p$, it is clear that
$$f_i(x)-f_i(x')\in {\cali}_n$$
where $f=(f_1,\ldots, f_p).$ Hence there exist $\alpha_{ij}(x,
x'),\;1\leq i\leq p,\;1\leq j\leq n,$ such that
$$f_i(x)-f_i(x')=\sum_{j=1}^{n}\alpha_{ij}(x,x')(x_j-x'_j).$$
If $f(x)=f(x')$ and $x\neq x'$, then clearly every $n\times n$
minor of the matrix $\alpha=(\alpha_{ij})$ must vanish at $(x,
x').$ We denote by ${\calr}_n(\alpha)$ the ideal in
${\calo}_{\C^{2n}}$ generated by the $n\times n$ minors of
$\alpha.$ Then we define the {\it double point ideal} as
$${\cali} ^2(f)=(f\times f)^{*}{\cali}_p+{\calr}_n(\alpha).$$

It is easy to verify that $V({\cali} ^2(f))=\tilde{D}^{2}(f)$ and
we call this complex space the {\it double point locus of} $f.$ At
a non-diagonal point $(x, x'),$ ${\cali} ^2(f)$ is generated by
the functions $f_i(x)-f_i(x').$ Moreover, the restriction of
${\cali} ^2(f)$ to the diagonal $\Delta_n$ is the ideal generated
by the $n\times n$ minors of the Jacobian matrix of $f$, so that
$\Delta_n\cap \tilde{D}^{2}(f)$ is just the singular locus of $f.$
The following property of the double point locus is a consequence
of \cite{Buchsbaum-Eisenbud}.

\begin{lema}\label{C-M}
The codimension of $\tilde{D}^{2}(f)$ is less than or equal to
$p$. Moreover, if the codimension is $p$, then $\tilde{D}^{2}(f)$
is Cohen-Macaulay.
\end{lema}

\begin{propo}\label{reduced}
Let $f:(\mathbb \C^{2},0)\rightarrow (\mathbb \C^{3},0)$ be a
finitely determined map germ. Then, $\tilde{D}^{2}(f)$ is a
reduced curve.
\end{propo}

\begin{proof}
By Lemma \ref{C-M} $\tilde{D}^{2}(f)$ is Cohen-Macaulay, hence it
is of pure dimension $1$ and satisfies Serre's conditions $(S_1)$,
that is all associated primes of ${\mathcal O}_{\tilde{D}^{2}(f)}$
are minimal (see for example \cite[p. 183]{Matsumura}). On the
other hand, by the Mather-Gaffney criterion for finite
determinacy, there is a representative $f:U\to \C^3$ of $f$
defined on some open neighborhood $U$ of $0$ in $\C^2$ such that
$f|_{U\setminus \{0\}}$ is stable and $f^{-1}(0)=\{0\}.$ Thus,
$\tilde{D}^{2}(f)\setminus \{0\}$ is a smooth curve in $U$ and
therefore ${\mathcal O}_{\tilde{D}^{2}(f)}$ satisfies Serre's
conditions $(R_0)$, that is the localization of ${\mathcal
O}_{\tilde{D}^{2}(f)}$ at every minimal prime is regular. Hence,
since ${\mathcal O}_{\tilde{D}^{2}(f)}$ satisfies both Serre's
conditions $(R_0)$ and $(S_1)$, $\tilde{D}^{2}(f)$ is reduced (see
for example \cite[p. 183]{Matsumura}).
\end{proof}

Assume that $G$ is a finite group which acts linearly on $\C^N$.
This action induces an analytic structure on the quotient $\C^N/G$
so that the local ring at a point $z\in\C^N$ is given by
$$\mathcal O_{N,z}^{G}=\{h\in\mathcal O_{N,z}: gh=h, \forall g\in G\}.$$
Assume now that $I\subset\mathcal O_{N,z}$ is a $G$-invariant
ideal. Then $G$ acts also on the germ of analytic set
$X=V(I)\subset (\C^N,z)$ and gives again an analytic structure on
$X/G$ with local ring $$\mathcal O_{X}^{G}=\{h\in\mathcal O_{X}:
gh=h, \forall g\in G\},$$ where $\mathcal O_X=\mathcal O_{N,z}/I$,
in such a way that $X/G$ embeds naturally in $(\C^N/G,z)$. If $I$
is generated by $G$-invariant functions $a_1,\dots,a_r\in\mathcal
O_{N,z}$, then $$\mathcal O_{X}^{G}\equiv \mathcal
O_{N,z}^{G}/I^G,$$ where $I^G$ is the ideal in $\mathcal
O_{N,z}^{G}$ generated by the same functions $a_1,\dots,a_r$.
Since $\mathcal O_{X}^{G}$ is in fact a subring of $\mathcal
O_{X}$, we have that if $X$ is reduced, then $X/G$ is also
reduced.

In our case, if $f$ is a holomorphic map or map germ from $\C^2$
to $\C^{3}$, then the double point ideal $\mathcal I^2(f)$ is
$S_2$-invariant, where we consider the action of the group $S_2$
on $\C^2\times \C^2$ given by $\tau(x,x')=(x',x)$. In this way, we
can define the quotient complex space or complex space germ
$\tilde{D}^2(f)/S_2$. It is a well known fact that $\C^2\times
\C^2/S_2$ is isomorphic to $\C^2$ times a quadratic cone in
$\C^3$. In particular, $\tilde{D}^2(f)/S_2$ embeds in $\C^5.$

If $f:U\subset\C^2\to\C^3$ is stable, then $\tilde{D}^2(f)$ is a
smooth curve and the quotient map $p:\tilde{D}^2(f)\to
\tilde{D}^2(f)/S_2$ is a 2-fold branched covering. For a finitely
determined map germ $f:(\C^2,0)\to(\C^3,0)$, $\tilde{D}^2(f)/S_2$
is the germ of a reduced curve and $\pi:\tilde{D}^2(f)\to
\tilde{D}^2(f)/S_2$ is a finite map germ, which is generically
2-to-1.

To complete the setup, let $D^2(f)=p_{1}(\tilde{D}^2(f))\subset
\C^2$ and $f(D^2(f))\subset \C^3,$ where $p_{1}:\C^2\times \C^2\to
\C^2$ is the projection on the first factor. For a finitely
determined map germ $f:(\C^2,0)\to(\C^3,0),$ $p_{1}$ is finite
when restricted to $\tilde D^2(f)$, in fact it is 1-to-1, and thus
$D^2(f)$ is the germ of a one-dimensional analytic set.
Analogously, $f$ is also finite when restricted to $D^2(f)$
(although in this case the map is 2-to-1 except at $0$), and thus
$f(D^2(f))$ is also a germ of a one-dimensional analytic set. In
both cases we consider the reduced analytic structure, so that
they become germs of reduced curves.

We have the following commutative diagram
\begin{equation}\label{CD}
\begin{CD}
\tilde{D}^2(f) @>>>   \tilde{D}^2(f)/S_2\\
@VVV @VVV\\
D^2(f) @>>> f(D^2(f)),
\end{CD}
\end{equation}
where the columns are 1-to-1.

Therefore, there are $3$ double-point sets associated to the
source: $\tilde D^2(f),$ $ D^2(f)$ and $ \tilde D^2(f)/{S_2},$ and
each of them is a reduced curve with isolated singularity. The
Milnor number of each of these sets, as defined in
\cite{Greuel-Buchweitz}, is an analytic invariant of the
singularity.

\begin{teore}[(Mond and Marar \cite{Mond-Marar})]\label{image-milnor}
Let  $f:(\C^2,0)\to (\C^3,0)$  be a finitely determined map germ,
$f_s$ a stable perturbation of $f$ and $X_s=f_s(\C^2)$, then the
following results hold.
 \begin{enumerate}
  \item $\chi(X_s)= C(f)+T(f)+\mu({\tilde D^2(f)}/{S_2})$,
  \item $X_s$ is simply connected and has the homotopy type of a
    wedge of $2$-spheres. The number of spheres in the
    wedge (known as the {\em{image Milnor number}}) is $$\mu_\Delta(f)= C(f)-1+T(f)+\mu( {\tilde D^2(f)}/{S_2}).$$
  \end{enumerate}
\end{teore}

\section{Families of finitely determined map germs }

\begin{defi}
Let $f:(\C^2,0)\to (\C^3,0)$ be a finitely determined map germ. A
\emph{one-parameter unfolding of} $f$ is a map germ
$F:(\C^2\times\C,0) \to (\C^3\times\C,0)$ of the form
$F(x,y,t)=(f_t(x,y),t)$ such that $f_t(0)=0, \, f_0=f$. We say
that an unfolding $F$ is a \emph{stabilization of} $f$ if there is
a representative $F:U\times T\to \C^3\times T,$ where $T$ and $U$
are open neighborhoods of $0$ in $\C$ and $\C^2$ respectively,
such that $f_t:U\to\C^3$ is stable for all $t\in T\setminus
\{0\}$.
\end{defi}

Since we are in the range of the nice dimensions in the sense of
Mather, it is well known that a stabilization of a finitely
determined map germ always exists.

Given an unfolding $F$ of $f$, we can also define the double point
locus $\tilde{D}^2(F)$ which is considered in $(\C^2\times
\C^2\times \C,0)$ instead of $(\C^3\times \C^3,0)$ and the other
set germs $\tilde{D}^2(F)/S_2$ in $(\C^2\times \C^2\times
\C/S_2,0)$, $D^2(F)$ in $(\C^2\times \C,0)$ and $F(D^2(F))$ in
$(\C^3\times \C,0)$.

\begin{lema}(\cite[Lemma 1.1]{Briancon-Galligo-Granger})\label{Briancon}
Let $X$ be a one parameter deformation of a reduced space curve
$X_0.$ Then, the following conditions are equivalent:

\begin{enumerate}

\item $X$ is a flat deformation of $X_0;$

\item $X$ is Cohen-Macaulay;

\item $X$ is of pure dimension $2.$

\end{enumerate}

\end{lema}

\begin{propo}\label{smooth} Let $F:(\C^2\times\C,0 )\to (\C^3\times\C,0)$
be an unfolding of a finitely determined map germ
$f:(\C^2,0)\to(\C^3,0).$  Then, the projections
$\tilde{D}^2(F)\to(\C,0)$, ${D}^2(F)\to(\C,0)$ and
$\tilde{D}^2(F)/S_2\to(\C,0)$ are flat deformations of
$\tilde{D}^2(f)$, ${D}^2(f)$ and $\tilde{D}^2(f)/S_2$
respectively. Moreover, if $F$ is a stabilization, then the
projection $\pi:\tilde{D}^2(F)\to(\C,0)$ is a smoothing.
\end{propo}

\begin{proof}
Let $\pi:\tilde{D}^2(F)\to(\C,0)$ be the projection. Note that the
second part of the proposition follows from the first one, since
for each $t\ne0$, $\tilde{D}^2(f_t)=\pi^{-1}(t)$ which is smooth
if $f_t$ is stable.

Let us show that $\pi:\tilde{D}^2(F)\to(\C,0)$ is a flat
deformation of $\tilde{D}^2(f)$. Note that $\tilde{D}^2(F)$ is
Cohen-Macaulay of pure dimension 2 by Lemma \ref{C-M}. We also
know by Proposition \ref{reduced} that $\tilde{D}^2(f)$ is reduced
and has pure dimension 1. Thus, $\tilde{D}^2(F)$ is a
Cohen-Macaulay $1$-parameter deformation of the reduced space
curve $\tilde{D}^2(f)$. Hence the result follows by Lemma
\ref{Briancon}.

The proof for $\tilde{D}^2(F)/S_2$ is analogous. Suppose that
$\mathcal I^2(F)^{S_2}$ is generated by $S_2$-invariant functions
$H_1,\dots,H_r\in\mathcal O_5$. Then $h_1,\dots,h_r\in\mathcal
O_4$ are also $S_2$-invariant, where $h_i(x,x')=H_i(x,x',0)$, and
generate $\mathcal I^2(f)^{S_2}$. Hence, $\mathcal I^2(F)^{S_2}$
restricted to $t=0$ is equal to $\mathcal I^2(f)^{S_2}$, that is
$\tilde{D}^2(F)/S_2$ is a $1$-parameter deformation of the space
curve $\tilde{D}^2(f)/S_2$. On the other hand, since
$\tilde{D}^2(F)$ is Cohen-Macaulay by Lemma \ref{C-M} we also have
that $\tilde{D}^2(F)/S_2$ is Cohen-Macaulay (see for example
\cite[Corollary 6.4.6]{Bruns-Herzog}). Thus, $\tilde{D}^2(F)/S_2$
is a Cohen-Macaulay $1$-parameter deformation of the reduced space
curve $\tilde{D}^2(f)/S_2$. Hence the result follows by Lemma
\ref{Briancon}.

Finally, the flatness of $D^2(F)$ follows from the fact that
$D^2(F)$ is a hypersurface in $\C^3$, hence it is Cohen-Macaulay.
Therefore, $\tilde{D}^2(F)$ is a Cohen-Macaulay $1$-parameter
deformation of the reduced plane curve $\tilde{D}^2(f)$. Hence the
result follows by Lemma \ref{Briancon}.
\end{proof}

\begin{defi} Let $F:(\C^2\times\C,0)\to(\C^3\times\C,0)$ be an unfolding of a finitely
determined map germ $f:(\C^2,0)\to(\C^3,0)$. We will say that $F$
is $\mu$-constant if $\mu(D^2(f_t))$ is constant along $T$, for
$t\in T$, where $T$ is a small neighborhood of $0$ in $\C.$
\end{defi}

\begin{defi} Let $F:(\C^2\times\C,0)\to(\C^3\times\C,0)$ be an unfolding of a finitely
determined map germ $f:(\C^2,0)\to(\C^3,0)$. We say that $F$ is
\emph{topologically trivial} if there are homeomorphism map germs:
\begin{align*}&\Phi:(\C^2\times \C,(0,0))\to (\C^2\times \C,(0,0)), &\Phi(x, t)=(\phi_t(x), t),\, \phi_0(x)=x,\,\phi_t(0)=0,\\
&\Psi:(\C^3\times \C,(0,0))\to (\C^3\times \C,(0,0)), &\Psi(y,
t)=(\psi_t(y), t),\, \psi_0(y)=y,\, \psi_t(0)=0,
\end{align*}
such that $F=\Psi \circ G\circ \Phi^{-1}$, where $G(x, t)=(f(x),
t)$ is the trivial unfolding of $f$.
\end{defi}

The following definitions were given by T. Gaffney in
\cite{Gaffney} for finitely determined map germs in dimensions
$(n,p).$ We restrict ourselves to the case $(n,p)=(2,3).$

\begin{defi}
  We say that $F$ is a \emph{good unfolding} of $f$ if there exist
  neighborhoods $U$ of $0$ in $\C^2\times\C$ and $W$ of $0$ in
  $\C^3\times\C,$ such that the following hold.
  \begin{enumerate}
  \item [(i)] $F^{-1}(\{0\}\times T)=\{0\}\times T,$ that is, $F$ maps
    $U\setminus(\{0\}\times T)$ into $W\setminus (\{0\}\times T)$.
  \item [(ii)] For all $(z_0,t_0)\in W\setminus (\{0\}\times T)$ the map
      $f_{t_0}:(\C^2,S)\to(\C^3,0)$ is stable, where
      $S=F^{-1}(z_0,t_0)\cap\Sigma F \cap U$ (which is a finite set) and
        $\Sigma F$ denotes the singular set of the unfolding $F$).
  \end{enumerate}
\end{defi}

\begin{defi}
  If there exists a curve $\alpha(s)=(x(s),y(s),t(s))$ in
  $(\C^2\times\C,0)$, containing $0$ in its closure, such that
  $(x(s),y(s))$ is a cross-cap point of $f_{t(s)}$, then we say that {\emph{$F$ has
   coalescing of cross-cap singularities}}.
\end{defi}

We can make an analogous definition for coalescing of triple
points. The following proposition is a particular case of
Proposition 3.6 in Gaffney's paper.

\begin{propo}
  Let $f:(\C^2,0)\to (\C^3,0)$  be a finitely
determined map germ and $F:(\C^2\times\C,0) \to (\C^3\times\C,0)$
be a good unfolding of $f.$ Then, the following conditions are
equivalent:
\begin{enumerate}
\item [(i)]$C(f_t)$ and $T(f_t)$ are constant. \item [(ii)] $F$
has no coalescing of cross-caps or triple points.
\end{enumerate}
\end{propo}

\begin{defi}
  A one-parameter unfolding $F$ of $f$ is an \emph{excellent
  unfolding} if it is good and the $0$-stable invariants $C(f_t)$ and $T(f_t)$ remain constant.
\end{defi}

Excellent unfoldings have a natural stratification. In the source
there are the following strata:
\[
\left\{ \C ^2 \setminus D^2(F), D^2(F)\setminus(\{0\}\times T),
\{0\}\times T \right\} .
\]
In the target, the strata are:
\[
\left\{ \C ^3 \setminus F(\C^2\times\C),
F(\C^2\times\C)\setminus\overline{F(D^2(F))}  ,
F(D^2(F))\setminus(\{0\}\times T), \{0\}\times T \right\} .
\]
Notice that $F$ preserves the stratification, that is, $F$ sends a
stratum into a stratum.

To finish the  section, we recall the main result in T. Gaffney's
paper \cite{Gaffney}, which gives, for any pair of dimensions
$(n,p),$ the topological triviality of excellent unfoldings of
finitely determined map germs for which the polar invariants
associated to the stable types in source and target remain
constant.

\begin{teore}[(Theorem 7.1, \cite{Gaffney})]
  Let $F:(\C^n\times\C,0) \to (\C^p\times\C,0)$ be an excellent
  unfolding. Let us assume that all polar invariants
of all stable types which occur in the stratification associated
to $F$ are constant at the origin. Then, the unfolding is
topologically trivial.
\end{teore}

For finitely determined map germs in the nice dimensions or their
boundary, Gaffney also proved (\cite{Gaffney}, Theorem 7.3) that
an excellent unfolding $F$ is Whitney equisingular if, and only
if, all polar invariants of all stable types of $F$ are constant.

\section{Properties of $\mu$-constant unfoldings}


\begin{teore}\label{teorema1}
  Let  $F:(\C^2\times\C,0) \to (\C^3\times\C,0)$ be a topologically
  trivial family of finitely determined map germs. Then,
  \begin{enumerate}
  \item [(a)] the Milnor numbers  $\mu(D^2(f_t)),\,\mu(\tilde D^2(f_t)),$ and $\mu({\tilde D^2(f_t)}/{S_2})$ are
  constant;
  \item [(b)] the $0$-stable invariants $C(f_t)$ and $T(f_t)$ are
  constant.
  \end{enumerate}
  In particular, $F$ is a $\mu$-constant unfolding of $f.$
\end{teore}

\begin{proof}
The following formulas were proved by Marar and Mond in
\cite{Mond-Marar}, Theorem 3.4:
\[
\mu(D^2(f_t))=\mu({\tilde D^2(f_t)})+6T(f_t)=2\mu({\tilde
D^2(f_t)}/{S_2})+C(f_t)+6T(f_t)-1.
\]
From these it suffices to prove that
$\mu(\tilde{D}^{2}(f_t)),\;\mu(\tilde{D}^{2}(f_t)/S_2)$ and
$\mu({D}^{2}(f_t))$ are constant. (We also use the fact that all
of the above invariants are upper semi-continuous.)

Since $F:(\C^2\times \C,0)\rightarrow (\C^3\times \C,0)$ defined
by $F(x,y,t)=(f_t(x,y),t)$ is topologically trivial,  there are
homeomorphisms $\Phi:(\C^2\times \C,0)\rightarrow (\C^2\times
\C,0)$ and $\Psi:(\C^3\times \C,0)\rightarrow (\C^3\times \C,0),$
as in the Definition 4.3, such that $ \Psi\circ F \circ \Phi^{-1}=
G, $
 where $G(x,y,t)=(f(x,y),t)$ is the trivial unfolding of
$f.$

It is easy to prove that these homeomorphisms give rise to
homeomorphisms (respecting the $S_2$-action) on the double point
locus:
 \begin{enumerate}
\item $\tilde{\Phi}^{2}:\tilde{D}^{2}(F)\rightarrow
\tilde{D}^{2}(G)$, \item
$\overline{\Phi}^{2}:\tilde{D}^{2}(F)/S_2\rightarrow
\tilde{D}^{2}(G)/S_2$, \item ${\Phi}^{2}:{D}^{2}(F)\rightarrow
{D}^{2}(G)$,
 \end{enumerate}
which obviously commute with the projections to $\C.$

 Clearly, $\tilde{D}^{2}(G)$ is homeomorphic to $\tilde{D}^{2}(f)\times
 T$, $\tilde{D}^{2}(G)/S_2$ to $\tilde{D}^{2}(f)/S_2\times
 T$ and ${D}^{2}(G)$ to ${D}^{2}(f)\times
 T$.

Hence, the above morphisms are topological trivializations of flat
families of reduced curves with isolated singularities. Therefore,
by Theorem 5.2.2 in \cite{Greuel-Buchweitz} we have that
$\mu({D}^{2}(f_t)),\,$ $\mu(\tilde{D}^{2}(f_t))$ and
$\mu(\tilde{D}^{2}(f_t)/S_2)$
 are constant.
\end{proof}

Next we will prove that $\mu$-constant unfoldings are excellent.

\begin{teore}\label{teorema2}
Let  $F:(\C^2\times\C,0) \to (\C^3\times\C,0)$ be an unfolding of
a finitely determined map germ $f.$ Then, the following statements
are equivalent:
\begin{enumerate}
\item\label{a} $\mu(D^2(f_t))$ is constant (i.e., $F$ is a
$\mu$-constant unfolding); \item\label{b}$\mu(\tilde D^2(f_t))$
and  $T(f_t)$ are constant; \item $\mu({\tilde D^2(f_t)}/{S_2}),$
$C(f_t)$ and  $T(f_t)$ are constant; \item The image Milnor number
$\mu_\Delta$  is constant.
\end{enumerate}
Moreover, any of these conditions implies that $F$ is an excellent
unfolding.
\end{teore}

\begin{proof}
The equivalence of the four conditions above comes from Theorem
\ref{image-milnor} and from the formulas from Marar and Mond
quoted in the proof of the previous theorem.

Let us assume now that $\mu(D^2(f_t))$ is constant in the
deformation. To prove that $F$ is excellent, we proceed as in
\cite{Gaffney}, Theorem 8.7.

To verify the first condition of goodness, we assume by
contradiction that,  in any neighborhood of $0\times 0$ in
$\C^2\times \C,$ we have $F^{-1}(\{0\}\times T)\neq \{0\}\times
T.$  If the points in $F^{-1}(\{0\}\times T)\setminus (\{0\}\times
T)$ lie in the singular set $\Sigma(F)$, of $F,$ then $C(f_t)$
must change at the origin, so we can assume they lie in
$F^{-1}(t)\setminus(\{0\}\times T \cup \Sigma (F)).$

Consider the intersection of $f_t(\C_x^2)$ with  $f_t(\C_0^2)$
where $(x,t)\in F^{-1}(\{0\}\times T)\setminus (\{0\}\times T)$.
If the intersection lies in $\Sigma(f_t(\C_0^2))$ then $T(f_t)$ is
at least one dimensional. Hence $f$ would not be finitely
determined if this holds for all $t$ close to $0.$ If
$f_t(\C_x^2)$ meets  $ \Sigma (f_t(\C_0^2))$ properly, then
$D^2(f_t)_x$ must have a singularity, hence $\mu(D^2(f_t))$ must
jump at $0.$

Considering the second condition for $F$ to be good, we suppose it
fails, so there exists an arc of points $(z(t),t)$ in
$\C^3\times\C,$ with $(0,0)$ in its closure such that $f_t$ is not
a stable map germ on $f_t^{-1}(z(t)).$ If  $f_t^{-1}(z(t))$
consists of a single point or three or more points, then
$z(t)\in\overline{f_t(\Sigma(f_t))\cup
  f_t(T(f_t))},$ so either $C(f_t)$  or $T(f_t)$ jumps at the origin.
The only possibility not so eliminated is that  $f_t^{-1}(z(t))$
consists of two points, say $x_1$ and $x_2$ with $x_1\neq x_2$, at
which $f_t$ is an immersion and $f(x_1)=f(x_2)$.

We now show that $x_1$ and $x_2$ are singular points of
$D^2(f_t).$ The easiest way  to see this is by picking disjoint
neighborhoods $U_i$ of $x_i$ and choosing coordinates centered at
$x_i$ so that $f_t|_{U_1}=(x_1,y_1,0)$ and
$f_t|_{U_2}=(x_2,y_2,f_3(x_2,y_2))$ with
$\mathrm{grad}\,f_3(0,0)=0.$ Then the equations for $\tilde
D(f_t)$ are
$$
\left\{
\begin{array}{rcl}
  x_1-x_2&=&0,\\
y_1-y_2&=&0,\\
f_3&=&0.
\end{array}
\right.
$$
So $\tilde D^2(f_t)$ has a singularity at $(0,0)\times (0,0)$ and
$\tilde D^2(f_t)$ sits in the diagonal of $U_1\times U_2.$ Since
the projection onto either factor is an isomorphism when
restricted to the diagonal, $D^2(f_t)$ in $\C^2$ has a singularity
at $x_i.$ Hence the constancy of $\mu(D^2(f_t))$ implies that such
an arc can not exist in this case either.
\end{proof}

\section{The main result}

\begin{lema}\label{Whitney}
Let  $F:(\C^2\times\C,0) \to (\C^3\times\C,0)$ be a $\mu$-constant
unfolding of a finitely determined map germ $f.$ Then, the
families of curves $D^{2}(F),\,\tilde{D}^{2}(F),$ and
$\tilde{D}^{2}(F)/S_2$ are Whitney equisingular along $\{0\}\times
T.$
\end{lema}

\begin{proof}
Since $\mu(D^{2}(f_{t}))$ is constant and $D^{2}(F)$ is a flat
family of reduced plane curves we have that  $(D^{2}(F)\setminus
(\{0\}\times T),\{0\}\times T)$ is Whitney regular. In particular,
the $D^{2}(f_{t})$ have constant multiplicity at $0$ (see
\cite{Greuel-Buchweitz}). On the other hand, by Theorem
\ref{teorema2} we have that the flat family $\tilde{D}^{2}(F)$ of
reduced curves have constant Milnor numbers and, since $p_1:
\C^2\times \C^2\times \C \to \C^2\times \C$ is a submersion which
restricts to a 1-to-1 map  $p_1:\tilde{D}^{2}(F)\to D^{2}(F),$ we
also have that the family $\tilde{D}^{2}(F)$ has constant
multiplicity at $0$. Hence, by a result of Buchweitz and Greuel
\cite{Greuel-Buchweitz} $(\tilde{D}^{2}(F)\setminus ( \{0\}\times
T),\{0\}\times T)$ is Whitney regular. Analogously,
$\tilde{D}^{2}(F)/S_2$ is also a flat family of reduced space
curves which has constant Milnor number and constant multiplicity,
because $p:\tilde{D}^{2}(F)\to \tilde{D}^{2}(F)/S_2$ is 2-to-1 and
the source has constant multiplicity. Hence,
$(\tilde{D}^{2}(F)/S_2\setminus (\{0\}\times T), \{0\}\times T)$
is Whitney regular.

\end{proof}

\begin{teore}\label{topological1}
Let $F:(\C^2\times\C,0) \to (\C^3\times\C,0)$ be a one-parameter
unfolding of a finitely determined map germ $f.$ The following
statements are equivalent:

\begin{enumerate}
\item $\mu(D^{2}(f_{t}))$ is constant for $t\in T,$ where $T$ is a
small neighborhood of $0$ in $\C;$ \item $F$ is topologically
trivial.
\end{enumerate}
\end{teore}
\begin{proof}
The implication  $(ii) \Rightarrow (i)$ follows from Theorem
\ref{teorema1}.

The idea of the proof of the converse statement is to construct
integrable vector fields $\xi$ in $\C^2\times \C$ and $\eta$ in
$\C^3\times \C,$ tangent to the strata of the stratifications in
source and target respectively,  such that
$$dF(\xi)=\eta \circ F.$$
In this way, the integral curves of $\xi$ and $\eta$ will define
the families of homeomorphisms $h_t:\C^2\to \C^2, \, k_t:\C^3 \to
\C^3$ such that
$$k_t \circ f_t \circ h_t^{-1}=f.$$

 Since ${D}^{2}(F)$ is Whitney equisingular along $\{0\}\times T$ as a
family of curves in $\C^2 \times \C,$ it follows from the First
Isotopy Lemma \cite{gwpl} that the vector field
${V}_{0}=\frac{\partial}{\partial t}$ on $\{0\}\times T$ lifts to
an integrable stratified vector field ${V}$ in a neighborhood $U$
of $0$ in $\C^2\times \C.$ The restriction of $V$ to each stratum
is smooth and tangent to the stratum, $d\pi_{T}(V)=V_{0}$ and
$d\rho(V)=0,$ where $\pi_{T}: \C^2\times \C \to \C$ is the
projection on the $t$-axis and $\rho$ is a control function,
$\rho: \C^2\times \C \to \C,$ $\rho \geq 0,\,\rho^{-1}(0)=
\{0\}\times T $ (see \cite{gwpl} for more details).

Notice that $D^2(F)$ has an analytic $S_2$-action given by $\sigma
(x,y,t)=(x',y',t),$ where $(x',y',t)$ is the unique point of
$D^2(F)$ such that $F(x,y,t)=F(x',y',t).$  We can average the
vector field $V$ to obtain a new controlled vector field on
$D^2(F),$
\[
W(x,y,t)= \frac{V(x,y,t)+ V(x',y',t)}{2},
\]
satisfying $W(x,y,t)=W(x',y',t),$ whenever $F(x,y,t)=F(x',y',t)$
and $d\pi_T(W)=\frac{\partial}{\partial t}.$ Notice that the
function $\rho^{*}(x,y,t)= \frac{\rho(x,y,t)+ \rho(x',y',t)}{2},$
is a control function for $W,$ that is $d\rho^{*}(W)=0.$

We can extend this vector field to a controlled vector field $\xi$
defined in a neighborhood of $0$ in $\C^2\times \C$, (for example
using the Kuo vector field, \cite{kuo}).

 The vector field $\eta$ given by $\eta=dF(\xi)$ is a well
  defined integrable vector field in the
image $F(D^2(F)).$ Furthermore, since $F$ is an isomorphism
outside the double point locus, $\eta$ can be extended to the
whole image $F(\C^2\times \C).$ Certainly we can extend this
vector field to the ambient space giving the desired topological
trivialization.
\end{proof}

\section{Bilipschitz triviality of $\mu$-constant unfoldings}

    A mapping $\phi : U\subset {K}^n\rightarrow {K}^p, \,\, K=\R  \, \text{or}\, \C,$ is called {\it Lipschitz} if
    there exists a constant $c >0$ such that:
\begin{equation*}
    \| \phi(x)-\phi(y) \| \leq c \| x-y \| \ \forall \ x,y
    \in U.
\end{equation*}
    When $n=p$ and $\phi$ has a Lipschitz inverse, we say that $\phi$
    is {\it bilipschitz}.

    Two map germs $f,g :({K}^n,0)\rightarrow ({K}^p,0)$ are called {\it
    bilipschitz equivalent} if there exists a bilipschitz map-germ $\phi
    :({K}^n,0)\rightarrow ({K }^n,0)$  such that
    $f=g\circ \phi .$

    A one-parameter unfolding $F:(\C^2\times\C,0) \to (\C^3\times\C,0)$
 of a finitely determined map germ $f: (\C^2,0) \to (\C^3,0)$
 is {\it bilipschitz trivial} if it is topologically trivial as in
 Definition 4.3 and the families of homeomorphisms in source and
 target are families of bilipschitz homeomorphisms.

When the bilipschitz trivialization is obtained by integrating
bilipschitz vector fields, we say that $F$ is {\it strongly
bilipschitz trivial}.

 Before proving the next theorem, we state a result of McShane
(\cite{MacShane}), which is a weak version of a Theorem of
Kirszbraun on the extension of Lipschitz functions:

\begin{teore}\label{Mac} (\cite{MacShane})
Let $X \subset \R^n$ be a metric subspace of the Euclidean space
and $f: X \to \R$ be a $c$-Lipschitz mapping, that is,
$||f(x)-f(y)|| \leq c ||x-y||.$ Then,
$$F(z)=sup_ {\{x \in X\}}(f(x)-c||x-z||),\, z\in \R^n$$
is a $c$-Lipschitz extension of $f.$

\end{teore}

\begin{rem}{\rm
If $f$ depends continuously on parameters $t=(t_1,\ldots,t_s),$
i.e., $f(x,t), \, x\in X$ is continuous in $(x,t)$ and is
$c$-Lipschitz in $x,$ with the constant $c$ not depending on $t,$
then it follows from the above result that there exists a
Lipschitz extension $F(x,t), \, x\in \R^n.$ Moreover, if $v$ is a
vector field on $X,$ depending continuously on parameters, we can
also apply Theorem \ref{Mac} to each coordinate function of $v$ to
obtain an extension $c\sqrt{n}$-Lipschitz $V$  to $\R^n.$}

\end{rem}

\begin{teore}\label{bilipschitz}
Let $F:(\C^2\times\C,0) \to (\C^3\times\C,0)$ be a one-parameter
unfolding of an  $\cala$-finitely determined map germ $f: (\C^2,0)
\to (\C^3,0).$ Then, the following are equivalent:

\begin{enumerate}
\item $F$ is topologically trivial. \item $F$ is Whitney
equisingular. \item $F$ is bilipschitz trivial.

\end{enumerate}

\end{teore}

\begin{proof}

To get the result it is sufficient to prove that $(a)$ implies
$(c),$ since the conditions $(c) \Longrightarrow (b)
\Longrightarrow (a) $ hold in general.  We do this following the
same steps as in Theorem \ref{topological}, constructing Lipschitz
vector fields $\xi$ in $\C^2\times \C,$ and $\eta$ in $\C^3\times
\C$ such that $dF(\xi)=\eta \circ F$.

As the set  $D^2(F)$ is a family of reduced plane curves,  then
the following are equivalent (cf.\ \cite{Greuel-Buchweitz,Za}):

\begin{itemize}

\item [(1)] $\mu(D^2(f_t))$ is constant;

\item [(2)] $D^2(F)$ is topologically trivial;

\item [(3)] $D^2(F)$ is Whitney equisingular;

\item [(4)] $D^2(F)$ is strongly bilipschitz trivial.

\end{itemize}

As in Theorem \ref{topological}, there is a vector field $V$
defined in a neighborhood of the origin in $\C^2\times \C,$
$d\pi_T(V)=\frac{\partial}{\partial t}.$  Since $D^2(F)$ satisfies
$(4)$, we now assume that $V$ is a Lipschitz vector field.

We now lift $V$ to a controlled  vector field $\tilde V$ in $(\C^2
\times \C^2 \times \C,0),$ which is tangent to the strata of
$\tilde D^2(F).$ Clearly, $d\pi_T(\tilde
V)=\frac{\partial}{\partial t}.$

Since the set $\tilde D^2(F)$ is $S_2$-invariant, we can average
the vector field $\tilde V$  to obtain a new controlled vector
field
$$\tilde W(x,y,x',y',t)= \frac{\tilde V(x,y,x',y',t)+\tilde
V(x',y',x,y,t)}{2},$$ tangent to the strata of the stratification
in $(\C^2 \times \C^2 \times \C,0),$ which is invariant under the
$S_2$ action on $\tilde D^2(F).$

Let $p_i: (\C^2 \times \C^2 \times \C,0) \to (\C^2 \times \C,0),\,
i=1,2, $ be the canonical projections: $p_1(x,y,x',y',t)=(x,y,t)$
and $p_2(x,y,x',y',t)=(x',y',t).$  The restriction of $p_1$ to
$\tilde D^2(F)$ is a $1$-to-$1$ generic projection in the sense
that the limit of all secant lines in $\tilde D^2(F)$ does not
live in the kernel of $p_1.$  Then, it follows from \cite{T}, page
$354,$ (see also \cite{alex}) that $p_1^{-1}: D^2(F) \to \tilde
D^2(F)$ is a bilipschitz homeomorphism. Then, $p_2\circ p_1^{-1} :
D^2(F) \to D^2(F)$ is also a bilipschitz homeomorphism. Notice
that $p_2\circ {p_1}^{-1}$ is the map-germ $\sigma$ defined in the
proof of Theorem~\ref{topological} and which gives the $S_2$
action in $D^2(F).$

Then the vector field  $W$ defined in $D^2(F)$ by $W(x,y,t)=
\tilde W(x,y, \sigma(x,y,t),t),$ is a Lipschitz vector field in
$D^2(F).$ We now use McShane's result to obtain a Lipschitz vector
field $\xi$ in $\C^2 \times \C$ which extends $W.$ The rest of the
proof follows as in Theorem 6.3, noticing that $\eta$ is clearly a
Lipschitz vector field, and that all the necessary extensions can
be made Lipschitz, applying again McShane's result.

\end{proof}

\begin{rem}{\rm
It is clear that Theorem \ref{topological} could be obtained as a
corollary of Theorem \ref{bilipschitz}. However, we present
Theorem \ref{topological} as our main result, because it has a
more general  setting, while Theorem \ref{bilipschitz} express
properties which do not hold in other dimensions.}

\end{rem}

As a consequence of the above theorem, it follows that the rich
theory of invariants of plane curves translate into results for
singular surfaces parametrized by finitely determined map germs.
We restrict ourselves to the following  corollary, which follows
from the above result, and Gaffney's Theorem 7.1, in
\cite{Gaffney}.

\begin{coro}

Let $F$ be as before  and $m_{0}(f_t)$ be the multiplicity of
$f_t$ at zero. If $F$ is topologically trivial, then  $m_0(f_t)$
is constant.
\end{coro}

\vspace{0.5cm}

{\small
\noindent Roberto Callejas Bedregal\hspace{3.4cm} Kevin Houston\\
{\it roberto@mat.ufpb.br} \hspace{4.3cm} {\it k.houston@leeds.ac.uk}\\
Universidade Federal da Para\'\i ba-DM \hspace{1.44cm}
University of Leeds\\
Cidade Universit\'aria s/n \hspace{3.71cm}School of Mathematics\\
58.051-900 Jo\~ao Pessoa, PB, Brazil. \hspace{1.6cm} Leeds, LS2
9J, United Kingdom.}

\vspace{1cm} {\small \noindent Maria Aparecida Soares Ruas\\
{\it maasruas@icmc.usp.br}\\
Departamento de Matem\'atica,\\
Instituto de Ci\^encias Matem\'aticas e de Computa\c{c}{\~a}o\\
Universidade de S\~ao Paulo-S{\~a}o Carlos\\
Caixa Postal 668, 13560-970 S\~ao Carlos SP, Brazil.}

\end{document}